\documentclass[generic,preprint]{imsart}
\RequirePackage{natbib}

\startlocaldefs
\usepackage{amssymb,amsmath,amsfonts,amsthm,amscd,amstext}

\usepackage[utf8]{inputenc}
\usepackage{latexmac}
\usepackage{graphicx}

\usepackage{url}
\usepackage{color}

\usepackage[charter]{mathdesign}

\usepackage[a4paper,scale={0.72,0.74},marginratio={1:1},footskip=7mm,headsep=10mm]{geometry}

\numberwithin{equation}{section}

\newtheorem{theorem}{Theorem}[section]
\newtheorem{claim}[theorem]{Claim}
\newtheorem{lemma}[theorem]{Lemma}
\newtheorem{proposition}[theorem]{Proposition}
\newtheorem{corollary}[theorem]{Corollary}
\newtheorem{remark}[theorem]{Remark}
\newtheorem{definition}[theorem]{Definition}

\long\def\xcom#1{}

\DeclareMathSymbol{\leqslant}{\mathalpha}{AMSa}{"36} 
\DeclareMathSymbol{\geqslant}{\mathalpha}{AMSa}{"3E} 
\DeclareMathSymbol{\eset}{\mathalpha}{AMSb}{"3F}     

\renewcommand{\epsilon}{\varepsilon}
\renewcommand{\theta}{\vartheta}
\renewcommand{\phi}{\varphi}

 \newcommand{\be}[1]{\begin{equation}\label{#1}}
 \newcommand{\ee}{\end{equation}}

 \newcommand{\bl}[1]{\begin{lemma}\label{#1}}
 \newcommand{\el}{\end{lemma}}

 \newcommand{\br}[1]{\begin{remark}\label{#1}}
 \newcommand{\er}{\end{remark}}

 \newcommand{\bt}[1]{\begin{theorem}\label{#1}}
 \newcommand{\et}{\end{theorem}}

 \newcommand{\bd}[1]{\begin{definition}\label{#1}}
 \newcommand{\ed}{\end{definition}}

 \newcommand{\bcl}[1]{\begin{claim}\label{#1}}
 \newcommand{\ecl}{\end{claim}}

 \newcommand{\bp}[1]{\begin{proposition}\label{#1}}
 \newcommand{\ep}{\end{proposition}}

 \newcommand{\bc}[1]{\begin{corollary}\label{#1}}
 \newcommand{\ec}{\end{corollary}}

 \newcommand{\bpr}{\begin{proof}}
 \newcommand{\epr}{\end{proof}}

 \newcommand{\bi}{\begin{itemize}}
 \newcommand{\ei}{\end{itemize}}

\xcom{
\newcommand{\esp}[1]{\mathbb{E}\etc{#1}}
\newcommand{\ens}[1]{\left\{#1\right\}}
\newcommand{\crochet}[1]{\left<#1\right>}

\newcommand{\etc}[1]{\left [#1 \right ]}
\newcommand{\etp}[1]{\left (#1 \right )}

\newcommand{\unsur}[1]{\frac{1}{#1}}

\newcommand{\undemi}{\frac12}
\newcommand{\un}[1]{1_{\etp{#1}}}

}

\date{\today}

\endlocaldefs

\begin{document}

\begin{frontmatter}
\runtitle{DPRE : two points state space}
\title{Directed polymer in random environment and two points state space}

\begin{abstract}

We give an exact expression  for the partition function of a continuous
time DPRE on a two points state space.
\end{abstract}

\author{\fnms{Philippe}
  \snm{Carmona}\corref{}\ead[label=e1]{philippe.carmona@univ-nantes.fr}
\ead[label=u1,url]{http://www.math.sciences.univ-nantes.fr/~carmona/}}

\affiliation{Université de Nantes}
\address{Laboratoire de Math\'ematiques Jean Leray UMR 6629\\
Universit\'e de Nantes, 2 Rue de la Houssini\`ere\\
BP 92208, F-44322 Nantes Cedex 03, France\\ \printead{e1}
}

\runauthor{P. Carmona et al.}



\end{frontmatter}

\section{Introduction}

We are aware of just three\footnote{If some interested reader finds
  another example where an exact computation of the free energy
  occurs, we are more than willing to incorporate it in this list}  models of directed polymer in a random
environment for which the free energy at finite inverse
temperature is well known:
\begin{enumerate}
\item The discrete polymer in dimension $d=1$ with simple random walk
  paths and 
  a log-gamma
  environment and boundary conditions: see \cite{MR2917766}.
\item The continuous time directed polymer with Poisson process paths
  and a Brownian environment
   : see \cite{MR1865759}.
\item The continuum random polymer with brownian paths and a white
  noise environment : see \cite{MR2796514}.
\end{enumerate}

The origin of this note is to find the simplest possible model where
some direct computations can be performed. Of course, this model is
much less interesting than the three models above. Note also that we
are even unable to generalize it to a three points state space !

Let us consider the simplest model of continuous time directed polymer
in a random environment. Let $\omega=(\omega(t))_{t\ge 0}$ be the
continuous time Markov chain on $\ens{1,2}$ with generator:
\begin{equation}
  Lf(1)= f(2) -f(1)\,,\quad Lf(2) = f(1) -f(2)\,,
\end{equation}
that is the chain that spends an exponential time on $1$ (resp $2$)
and the jumps on $2$ (resp. $1$). We let $P_i$ denote the law of the
Markov chain starting from $i$, $E_i$ the associated expectation, and
$(W,\Wrond)$ the path space of piecewise constant cadlag functions
from $[0,+\infty[$ to $\ens{1,2}$. We set $P=P_1$ and $E=E_1$.

The random environment consists of two independent standard Brownian
motions $(B_i(t), t\ge 0)$ defined on another probability space
$(\Omega,\Frond,\PP)$.

For any $t>0$ the (random) \emph{polymer measure} $\mu_t$ is the probability defined on the path space $(W,\Wrond)$ by
\begin{equation*}
  \label{eq:1}
  \mu_t(d\omega)= \unsur{Z_t} e^{\beta H_t(\omega) - t \beta^2/2} \PP(d\omega)
\end{equation*}
where $\beta \ge 0$ is the inverse temperature, the Hamiltonian is
\begin{equation*}
  \label{eq:2}
  H_t(\omega) = \intot dB_{\omega(s)}(s)
\end{equation*}
and the \emph{partition function} is 
\begin{equation*}
  \label{eq:3}
  Z_t = Z_t(\beta) = E\etc{e^{\beta H_t(\omega) - t \beta^2/2}}\,,
\end{equation*}

It is well known, and for example established in \cite[Proposition
2.4]{MR2249669} that the \emph{point to point partition function}
\begin{equation}
  Z_t(x,y) := E_x\etc{e^{\beta H_t(\omega) - t \beta^2/2} \un{\omega(t)=y}}
\end{equation}
satisfy the \emph{Discrete Stochastic Heat Equation}:
\begin{equation}
  dZ_t(x,y) = LZ_t(x,.)(y)\, dt + \beta Z_t(x,y) dB_y(t)\,,
\end{equation}
in the Itô sense with inital conditions
$Z_0(x,y)=\delta_x(y)$ (from now on we shall fix as a starting point $x=1$). Furthermore, the free energy is well defined
\begin{equation}
  p(\beta) := \lim_{t\to +\infty} \unsur{t} \log Z_t(\beta)
  \qquad(\text{a.s. and in $L^1(\PP)$})\,,
\end{equation}
and is given by the limit (see \cite[Formula (15)]{MR2249669})
\begin{equation}
  p(\beta)= - \frac{\beta^2}{2} \lim_{t\to +\infty} \unsur{t} \int_0^t
  \esp{I_s}\, ds\,,
\end{equation}
with $I_t$ the\emph{ overlap}
$$I_t=\mu_t^{\otimes 2}\etp{\omega_1(t)=\omega_2(t)} =\unsur{Z_t^2}
\sum_{y=1}^2 Z_t(1,y)^2\,.$$

\begin{theorem}\label{thm:un}
  For this model of DPRE the mean ovelap converges:
  \begin{equation}
    \lim_{t\to +\infty} \esp{I_t} = \alpha_-(\beta)
  \end{equation}
with $\alpha_-(\beta) \in (0,1)$ the smallest root of the polynomial
$$ P_\beta(X) = 3 \beta^2 X^2 - (5\beta^2 + 4) X + 2 (1+\beta^2)\,.$$
Consequently the free energy is
$$ p(\beta) = - \frac{\beta^2}{2} \alpha_-(\beta)\,.$$
\end{theorem}

Observe that as $\beta\to 0^+$, $\alpha_-(\beta) \to \undemi$ as
expected.

\section{Proof of Theorem~\ref{thm:un}}

To simplify notations we let $X_i(t)=Z_t(1,i)$ and set $Z_0=1$ so that
we have $X_1(0)=1$, $X_2(0)=0$, $I_0=1$ and $(X_1,X_2)$ is solution of
the following \emph{simple} system of stochastic differential equations:

\begin{equation}\label{eq:simplesystem}
\left\{
\begin{split}      dX_1(t)& = (X_2 -X_1)\, dt + \beta X_1 dB_1(t)\\
dX_2(t) &= -(X_2 -X_1)\, dt + \beta X_2 dB_2(t)
\end{split}
\right.
\end{equation}

It is easy to check that $Z_t= X_1(t) + X_2(t)$ is a martingale 
$$ dZ_t = \beta (X_1 dB_1(t) + X_2 dB_2(t))$$
with quadratic variation
$$d\crochet{Z}_t =\beta^2(X_1^2 + X_2^2)\, dt =
  \beta^2 N_t \, dt = \beta^2 Z_t^2 I_t\, dt\,,
$$
where we have set 
$$N_t =X_1(t)^2 + X_2(t)^2\quad\text{ so that }\quad I_t = \frac{N_t}{Z_t^2}\,.$$

Without having to read \cite{MR2249669}, one can infer directly that
$$ \log Z_t = \int_0^t \frac{dZ_s}{Z_s} - \undemi \int_0^t
\frac{d\crochet{Z}_s}{Z_s^2} = \int_0^t \frac{dZ_s}{Z_s} - \frac{\beta^2}{2}
\int_0^t I_s\, ds$$
so that 
$$ \unsur{t} \esp{\log Z_t} = - \frac{\beta^2}{2} \unsur{t} \int_0^t
\esp{I_s}\, ds\,.$$

Let us do now some straightforward computations using Ito's formula
\begin{align*}
d X_1^2(t) &= ((\beta^2-2) X_1^2 + 2 X_1 X_2)\, dt + 2 \beta X_1^2 dB_1\\
dN_t &= (4 X_1X_2 + (\beta^2-2) N_t) dt + 2 \beta (X_1^2 dB_1 + X_ 2^2
dB_2) = (2 Z_t^2 + (\beta^2-4)N_t)\, dt +  2 \beta (X_1^2 dB_1 + X_ 2^2
dB_2)\\
 d\crochet{N,Z}_t &= 2 \beta^2(X_1^3 + X_2^3) dt =
\beta^2(3 Z_t N_t - Z_t^3)\, dt = \beta^2 Z_t^3(3 I_t-1)\, dt\,.
\end{align*}
In the last equation, we use the identity $2 (a^3 + b^3) = 3 (a+b)
(a^2 + b^2) - (a+b)^3$. It is the only place where we really use the
fact that the state space has only two points.

Let use the notation $U_t \sim V_t$ if $U_t -V_t$ is a martingale. We have:

\begin{align}
  dI_t &= \frac{dN_t}{Z_t^2} - \frac{2 N_t dZ_t}{Z_t^3}
  -2\frac{d\crochet{N,Z}_t}{Z_t^3} + 3 N_t
  \frac{d\crochet{Z}_t}{Z_t^4}\notag\\
&\sim ((\beta^2-4)I_t + 2 - 2 \beta^2(3 I_t -1) + 3 \beta^2 I_t^2)\,
dt
\notag\\
&\sim (3 \beta^2 I_t^2 - (5\beta^2+4) I_t + 2 (1+\beta^2))dt =
P(I_t)\, dt \label{eq:dit}
\end{align}
with $a=\frac{5\beta^2 + 4}{6 \beta^2} = \frac56 + \frac{2}{3\beta^2}$
and 
\begin{align*}
  P(X) = 3 \beta^2(X^2 - 2a X + a- \unsur{6}) = 3 \beta^2(X-\alpha_+)(X-\alpha_-)
\end{align*}
with $\alpha_\pm := a \pm \sqrt{a^2 -a + \unsur{6}}$. Since $P(1)=3
\beta^2(1-a -\unsur{6}) = -2$ we have  $\alpha_- < 1 < \alpha_+$.

We now take expectations in \eqref{eq:dit} and get

\begin{equation} \frac{d \esp{I_t}}{dt} = \esp{P(I_t)}\,.
\end{equation}
Since $I_0=1$, and $P(1)<0$, at least on a non empty  interval
$[0,\delta[$, the function  $u(t):=\esp{I_t}$ is non increasing
(something we expected since 
$0\le I_t \le 1$ donc $0\le u(t) \le 1$).

Assume that there exists  $t_0>0$ such that $u(t_0) <
\alpha_-$. Since $\esp{I_t^2}\ge \esp{I_t}^2$ we have

$$ u'(t) = \esp{P(I_t)} \ge P(u(t))\,.$$

Let $T=\sup\ens{t< t_0 : u(t)=\alpha_-}$. On $]T,t_0[$ we have $u'(t)\ge
P(u(t)) >0$ so $u$ strictly increases on $]T,t_0[$, which is absurd
since $u(T)=\alpha_-$ and $u(t_0)< \alpha_-$.

 Therefore, we have:
$$ \forall t>0, u(t) \ge \alpha_-\,.$$

Observe that:
\begin{align*}
  u'(t) &= \esp{P(I_t)} = \esp{P(I_t) -P(\alpha_-)} =
  \esp{3\beta^2(I_t^2-\alpha_-^2) - (5\beta^2+4) (I_t -\alpha_-)} \\
& = 3 \beta^2 \esp{(I_t + \alpha_-)(I_t -\alpha_-)} - (5\beta^2 + 4)
\esp{I_t- \alpha_-}\\
&\le \etp{3 \beta^2(1+\alpha_-) - (5\beta^2 + 4) }\esp{I_t- \alpha_-} =
-\lambda \esp{I_t- \alpha_-}\,.
\end{align*}
We shall check that $\lambda >0$, and this implies, by Gronwall's
Lemma, 
$$ (u(t) -\alpha_-) \le (u(0)-\alpha_-)e^{-\lambda t} \to 0
\quad\text{when $t\to +\infty$}\,.$$

It remains to prove that 
 $\lambda = -3 \beta^2(1+\alpha_-) + (5\beta^2 +
4) >0$ that is  $1+\alpha_- < 2a$, i.e. $\sqrt{a^2 -a +
  \unsur{6}} > 1-a$.
 
Either $\beta \le 2$ and then $a\ge 1$ and we are done, or
 $\beta >2$, $a<1$, and we have to show that   $a^2 -a +
\unsur{6} > (1-a)^2$ i.e. $a>5/6$ which is true.

\bibliographystyle{imsart-nameyear}
\bibliography{amscar}

\end{document}